\documentclass{conm-p-l}

\usepackage{amssymb}

\newtheorem{theorem}{Theorem}[section]
\newtheorem{lemma}[theorem]{Lemma}

\theoremstyle{definition}
\newtheorem{definition}[theorem]{Definition}
\newtheorem{example}[theorem]{Example}
\newtheorem{proposition}[theorem]{Proposition}

\theoremstyle{remark}
\newtheorem{remark}[theorem]{Remark}

\numberwithin{equation}{section}

\newcommand{\nc}[2]{\newcommand{#1}{#2}}
\newcommand{\rnc}[2]{\renewcommand{#1}{#2}}
\rnc{\theequation}{\thesection.\arabic{equation}}
\nc{\beq}{\begin{equation}}
\nc{\eeq}{\end{equation}}

\nc{\bpr}{\begin{proposition}}
\nc{\bth}{\begin{theorem}}
\nc{\ble}{\begin{lemma}}
\nc{\bco}{\begin{corollary}}
\nc{\bre}{\begin{remark}}
\nc{\bex}{\begin{example}}
\nc{\bde}{\begin{definition}}
\nc{\ede}{\end{definition}}
\nc{\epr}{\end{proposition}}
\nc{\ethe}{\end{theorem}}
\nc{\ele}{\end{lemma}}
\nc{\eco}{\end{corollary}}
\nc{\ere}{\end{remark}}
\nc{\eex}{\end{example}}
\nc{\bpf}{\begin{proof}}
\nc{\epf}{\end{proof}}

\nc{\mc}{\mathcal{C}}
\nc{\wmc}{\widetilde{\mathcal{C}}}
\nc{\ws}{\widetilde{S}}
\nc{\wdel}{\widetilde{\delta}}
\nc{\wsigma}{\widetilde{\sigma}}
\nc{\wtau}{\widetilde{\tau}}
\nc{\wvfy}{\widetilde{\varphi}}

\nc{\ra}{\rightarrow} \nc{\lra}{\longrightarrow}
\nc{\lla}{\longleftarrow}

\nc{\epn}{\varepsilon}
\nc{\si}{\psi}
\nc{\Si}{\Psi}
\nc{\Del}{\Delta}
\nc{\del}{\delta}
\nc{\ro}{\rho}
\nc{\fy}{\phi}
\nc{\Fy}{\Phi}
\nc{\tu}{\tau}
\nc{\vfy}{\varphi}
\nc{\sig}{\sigma}

\nc{\tta}{\theta}
\nc{\alf}{\alpha}
\nc{\bt}{\beta}
\nc{\mh}{\mathcal{H}}

\nc{\mbc}{\mathbb{C}}
\nc{\mbz}{\mathbb{Z}}
\nc{\mbr}{\mathbb{R}}
\nc{\mfg}{\mathfrak{g}}
\nc{\ug}{U(\mfg)}

\nc{\nn}{\nonumber}

\nc{\ot}{\otimes}

\newcommand{\lrbicross}{{\blacktriangleright\!\!\!\triangleleft}}

\nc{\dff}{Diff^+(\mbr^{1,1})}
\nc{\aff}{Aff(\mbr^{1,1})}
\nc{\rnd}{\partial}
\nc{\rr}{\mbr^{1,1}}

\begin{document}

\title{A Super Version of the Connes-Moscovici Hopf Algebra}

\author{Masoud Khalkhali}
\address{Department of Mathematics,  University of Western Ontario,
London, Ontario, Canada, N6A 5B7}
\email{masoud@uwo.ca}

\author{Arash Pourkia}
\address{Department of Mathematics,  University of Western Ontario,
London, Ontario, Canada, N6A 5B7}
\email{apourkia@uwo.ca}

\date{January 1, ... and, in revised form, June 22, ....}

\dedicatory{Dedicated, with much appreciation,  to Henri Moscovici.}


\begin{abstract}
We define a super version of the Connes-Moscovici Hopf algebra, $\mh_1$.
For that, we consider the supergroup $G^s=Diff^+(\mbr^{1,1})$ of
orientation preserving diffeomorphisms of the superline $\mbr^{1,1}$ and
define two (super) subgroups $G^s_1$ and $G^s_2$ of $G^s$ where $G^s_1$ is
the supergroup of affine transformations. The super Hopf algebra $\mh ^s_1$ is defined
as a certain \emph{bicrossproduct super Hopf algebra} of the super Hopf algebras attached to
$G^s_1$ and $G^s_2$. We also give an explicit description of $\mh ^s_1$ in terms of
generators and relations.
\end{abstract}

\maketitle

\section{Introduction}

In \cite{cm2}, Connes and Moscovici, among many other things,  defined a Hopf algebra $\mh(n)$, for any $n \geq 1$,
and computed the periodic Hopf cyclic cohomology of
 $\mh(n)$.
Our main focus in this paper is $\mh(1)$, to be denoted from now on  by $\mh_1$, and its
 super analogue. It is by now clear that the Connes-Moscovici Hopf algebra $\mathcal{H}_1$ is a fundamental object
 of noncommutative geometry.   An
 important feature of $\mathcal{H}_1$, and in fact its  {\it raison
 d'\^etre}, is that it acts as  quantum symmetries of  various algebras  of interest in noncommutative
 geometry, like the  algebra of leaves of codimension one foliations and  the algebra  of modular forms
 modulo the action of Hecke correspondences \cite{cm2, cm4, cm5, cm6}.

Our starting point,  in fact the motivation to develop the super analogue of $\mathcal{H}_1$,  was
to extend the results of \cite{cm6} to cover the Rankin-Cohen brackets on super modular forms and super pseudodifferential operators as they are
described in Section 7 of \cite{cmz}.   In  \cite{cm6} it is shown that  the Rankin-Cohen brackets on modular forms \cite{cmz} can
be derived via  the action of $\mathcal{H}_1$ on the modular  Hecke algebras. In fact, more generally, it is shown how to obtain  such brackets on
 any associative algebra endowed with an action of the Hopf algebra $\mathcal{H}_1$,
  such that the derivation corresponding to the Schwarzian derivative is
  inner. To carry out this program in the SUSY case, as a first step one
  needs a super analogue of $\mathcal{H}_1$.

 The Connes-Moscovici Hopf algebra $\mathcal{H}_1$
 is isomorphic to a certain {\it{bicrossproduct}} Hopf algebra
$F(G_2) \lrbicross U(\mfg_1)$  \cite{cm2,cm4,maj3,bh07,bh09}. The actions and coactions
involved in this  bicrossproduct can be derived and understood  by looking at
 the factorization of the group of orientation preserving diffeomorphisms
of the real line, $G=Diff^+(\mbr)$, into two subgroups $G_1$ and $G_2$. Here
$G_1$ is the group of affine transformations and $G_2$ is the subgroup of those
diffeomorphisms $\fy$ with $\fy(0)=0$ and $\dot{\fy}(0)=1$, where
$\dot{\fy}(x)= \frac{d}{dx}(\fy(x))$.

In this paper our goal is to define a super version of  $\mh_1$, which we
will denote by  $\mh ^s_1$. For that  we define
a super version of the group $G=Diff^+(\mbr)$, namely the supergroup
$G^s=\dff$ of orientation preserving diffeomorphisms of the superline $\mbr^{1,1}$.
We define two (super) subgroups $G^s_1$ and $G^s_2$
of $G^s$,  where $G^s_1$ is the group of affine transformations.
We show that the factorization $G^s=G^s_1G^s_2$ holds.
We will use this factorization to define a  bicrossproduct super Hopf algebra
$F(G^s_2) \lrbicross U(\mfg^s_1)$, analogous to the non-super case.
We will call this bicrossproduct super Hopf algebra
the super version of  $\mh_1$ and denote it by $\mh ^s_1$.

One difficulty in working with super Hopf algebras is that they are not `honest' Hopf algebras. In fact their multiplication map is not a morphism
of coalgebras. It is so only up to sign, and this issue of signs can be quite confusing and demands a lot of care.
Throughout this paper, for notations to be more consistent
with the non-super case, we use the following conventions. To denote the comultiplication
$\Del: B \ra B\ot B$ of a bialgebra $B$ we use  Sweedler's notation (summation understood)
$\Del^n(b)=b_{(1)}\ot b_{(2)}\ot \cdots \ot b_{(n+1)}$, for any $n \geq 1$. Also for a
coaction $\nabla : A \ra B\ot A$ of $B$ on $A$ we write  $\nabla (a)= a^{(1)} \ot  a^{(2)}$.

One of the recent important developments regarding the Connes-Moscovici Hopf algebras $\mh(n)$ is in \cite{bh09}, in which the authors generalize the Connes-Moscovici  work, \cite{cm2}, to assign a Hopf algebra to any infinite primitive Lie pseudogroup. They also introduce an elaborate machinery, based on the bicrossproduct realization of those Hopf algebras, to compute their periodic and non-periodic Hopf cyclic cohomology.

We would like to thank Bahram Rangipour for useful discussions and  for suggesting a  proof to show that $F(G^s_2)$ is a super Hopf algebra. We would also like to thank Arthur Greenspoon for carefully reading the manuscript and many useful suggestions.  
We are honored  to dedicate this paper,  with much appreciation,  to
 Henri Moscovici on the occasion of his 65th birthday.

\section{The Connes-Moscovici Hopf algebra $\mh_1$}
In this section we recall the definition of the standard (non-super) Connes-Moscovici Hopf algebra $\mh_1$ and its description in terms of a
 bicrossproduct Hopf algebra  \cite{cm2,cm4}. The following definition gives a description of $\mh_1$ by generators and relations.
\bde  \label{definofCMH1} \cite{cm2,cm4}
The Connes-Moscovici Hopf algebra  $\mh_1$ is generated by elements $X$, $Y$, $\del_n,\,\,n \geq 1$ with relations:
\[ [Y,X]=X,\,\,\,[X,\del_n]=\del_{n+1},\,\,\,[Y,\del_n]=n \del_n,\,\,\,[\del_m,\del_n]=0,\quad \forall m,n \nn\]
\[\Del(X)=X \ot1 +1\ot X+Y\ot \del_1 ,\,\,\,\Del(Y)=Y \ot1 +1\ot Y,\,\,\,\Del(\del_1)=\del_1 \ot1 +1\ot  \del_1 , \nn\]
\[\epn(X)=\epn(Y)=\epn(\del_n)=0,\,\,\forall n , \nn\]
\beq
S(X)=Y\del_1 - X,\,\,\,S(Y)=-Y, \,\,\,S(\del_1)=-\del_1 .
\eeq
\ede

\bre
The above definition follows the right-handed notation as in \cite{maj3}, in the sense that in the definition of $\Del(X)$, the term $\del_1 \ot Y$
in  \cite{cm2,cm4} is replaced by  $Y\ot \del_1$. Or, alternatively, if we denote the original Connes-Moscovici Hopf algebra by ${\mh_1}_{CM}$, then one can say we are working in ${\mh_1}_{CM}^{cop}$, \cite{bh09}.
\ere

\ble  \cite{maj3,maj}  \label{cptbltycondnonsup}
let $A$ and $H$ be two Hopf algebras such that $A$ is a left $H$-module algebra, and $H$ is  a right $A$-comodule coalgebra. Let furthermore these structures
 satisfy the follow compatibility conditions:
\[\Del(h \rhd a)= {h_{(1)}}^{(1)} \rhd a_{(1)} \ot {h_{(1)}}^{(2)} (h_{(2)} \rhd a_{(2)})    , \nn \]
\[\nabla_r(gh)={g_{(1)}}^{(1)} h^{(1)} \ot {g_{(1)}}^{(2)}(g_{(2)} \rhd h^{(2)}) , \nn \]
\[{h_{(2)}}^{(1)} \ot ({h_{(1)}}\rhd a){h_{(2)}}^{(2)} ={h_{(1)}}^{(1)} \ot {h_{(1)}}^{(2)} (h_{(2)} \rhd a)    , \nn \]
\[\epn(h\rhd a)=\epn(h) \epn(a), \quad \quad \nabla_r(1)=1\ot 1 , \nn\]
 for any $a,b$ in $A$ and $g,h$ in $H$, where we have denoted the actions by $h \rhd a$ and the coactions by
$\nabla_r (h)= h^{(1)} \ot  h^{(2)}$.
Then the vector space $A\ot H$ can be equipped with a Hopf algebra structure as follows:
\[(a\ot h)(b\ot g)= a(h_{(1)} \rhd b) \ot h_{(2)} g     , \nn \]
\[\Del(a \ot h)=a_{(1)} \ot {h_{(1)}}^{(1)} \ot a_{(2)}{h_{(1)}}^{(2)}\ot h_{(2)}     , \nn \]
\[\epn(a \ot h)= \epn(a)\epn(h)     , \nn \]
\[S(a \ot h)= (1\ot S(h^{(1)}))(S(a h^{(2)}) \ot 1)     , \nn \]
for any $a,b$ in $A$ and $g,h$ in $H$.
 It is called the left-right bicrossproduct Hopf algebra and is denoted by  $A \lrbicross H$.
\ele

in  \cite{cm2,cm4} a Hopf subalgebra  of  $\mh_1$ is defined as the unital commutative subalgebra of  $\mh_1$ generated by
$\{ \del_n, n \geq 1 \}$.
   This Hopf algebra, which we denote by  $F(G_2)$, is isomorphic to the  so-called {\it{comeasuring}} Hopf algebra
of the real line, generated by
$\{a_n, \quad n \geq 1 \}$ with $a_1=1$ and with the following relations \cite{fg,maj3}:
\[\Del(a_n)= \sum_{k=1}^n (\sum_{i_1 +\cdots +i_k =n} a_{i_1} \cdots a_{i_k}) \ot a_k , \nn \]
\[\epn(a_n)=\del_{n,1} , \label{fg2def} \nn \]
\[S(a_{n+1})= \sum_{(c_1, \dots,c_{n+1}) \in \Lambda }
          (-1)^{n-c_1} \frac{(2n-c_1)!c_1 ! \quad a_1^{c_1} a_2^{c_2} \cdots a_{n+1}^{c_{n+1}}}{(n+1)! \quad
c_1!c_2! \cdots c_{n+1}!}    , \nn \]
where
\[\Lambda = \{(c_1, \dots,c_{n+1}) \quad | \quad \sum_{j=1}^{n+1} c_j =n, \quad  \sum_{j=1}^{n+1} jc_j =2n  \} . \nn \]
If instead of generators $a_n$ we work with $n!a_n$ we get the relations for the so called {\it{Fa\`a di Bruno}} Hopf algebra which is isomorphic to  $F(G_2)$.
 One can also define another Hopf algebra, denoted by  $U(\mfg_1)$, as follows. Let $\mfg_1$ be the Lie algebra
generated by two elements $X$ and $Y$ as in  $\mh_1$
(i.e., $[Y,X]=X$), and let  $U(\mfg_1)$ denote the universal enveloping algebra of  $\mfg_1$.

\ble \label{fg2actandg1coact} \cite{maj3}
 $F(G_2)$ is a left $U(\mfg_1)$-module algebra via the actions
\[X\rhd a_n=(n+1)a_{n+1} -2a_2 a_n , \quad Y \rhd a_n=(n-1)a_n  , \nn \]
and  $U(\mfg_1)$ is a right $F(G_2)$-comodule coalgebra via the coactions
\[\nabla_r(X)=X\ot 1 +Y\ot 2a_2 , \quad \nabla_r(Y)=Y\ot 1  . \nn \]
The actions and coactions are compatible in the sense of Lemma \eqref{cptbltycondnonsup}.
\ele

\bth \cite{cm2,cm4,maj3, bh07,bh09}
 $\mh_1$ is isomorphic to the  bicrossproduct Hopf algebra $F(G_2) \lrbicross  U(\mfg_1)$.
\ethe

A good way to understand the actions, coactions and even the notations  introduced above  is to look at the \emph{factorization} of the group $Diff^+(\mbr)$ \cite{cm2,cm4}.
Let us recall that given any group $G$ with two subgroups $G_1$ and $G_2$, we say we have a \emph{group factorization} $G=G_1G_2$ if
 for any $g \in G$ there is a unique decomposition $g=ab$ where $a\in G_1$ and $b\in G_2$.
Given any group factorization $G=G_1G_2$, one has always a left action of $G_2$ on $G_1$ and a right action of
$G_1$ on $G_2$ defined in the following way. First one defines two maps $\pi_1 :G \ra G_1$,  $\pi_2 :G \ra G_2$ by $\pi_1(g)=a$ and $\pi_2(g)=b$, for any $g=ab$ in $G$, with $a\in G_1$ and $b\in G_2$.
Next, one can define the aforementioned actions by $g_2 \rhd g_1=\pi_1 (g_2g_1)$ and  $g_2 \lhd g_1=\pi_2 (g_2g_1)$ for any $g_1$ in $G_1$ and $g_2$ in $G_2$.

Let \[G=Diff^+(\mbr)=\{\Fy \in Diff(\mbr) \,| \, \dot{\Fy}(x) > 0,\, \forall x \in \mbr \} ,\nn\]
where $\dot{\Fy}(x)= \frac{d}{dx}(\Fy(x))$,
be the group of orientation preserving homeomorphisms of the real line and
\[G_1= \{\si=(a,b) \in G \, | \,  \si(x)=ax+b , \, a,b \in \mbr,\, a>0\}, \nn \]
be the affine subgroup of $G$. The following representation of $G_1$ as a subgroup of $GL(2)$ is very useful:
\beq\label{affingl2}
G_1=\left\{ (a,b)=\left( \begin{array}{cc}
a & b  \\
0 & 1  \\
\end{array} \right) \in GL(2) \,| \, a>0 \right\}.
\eeq
Let also
\[G_2=\{\fy \in G \, | \, \fy(0)=0,\, \dot{\fy}(0)=1 \}. \nn \]
The factorization $Diff^+(\mbr)=G_1G_2$ is as follows. For any $\Fy$ in $G$ we have $\Fy=\si \fy$ with $\si \in G_1$
and $\fy \in G_2$, where
\beq\label{facrofG}
\si=(\dot{\Fy}(0), \Fy(0))  , \quad \quad \fy(x)=\frac{\Fy(x)-\Fy(0)}{\dot{\Fy}(0)}, \quad \forall x \in \mbr .
\eeq

The actions and coactions in lemma \eqref{fg2actandg1coact} are induced from the factorization of the group $G=Diff^+(\mbr)$ as follows:
\[\fy \rhd \si = (a\dot{\fy}(b),\fy(b)) ,\quad \quad (\fy \lhd \si)(x)=\frac{\fy(ax +b) -\fy(b)}{a \dot{\fy}(b)} , \nn \]
for any $\si=(a,b)$ in $G_1$ and $\fy$ in $G_2$.

Next, using the matrix representation of $G_1$ in $GL(2)$ as in \eqref{affingl2},
we get a basis for its Lie algebra $\mfg_1=Lie(G_1)$. It turns out that
$\mfg_1$ is generated by two elements
\[X= \left( \begin{array}{cc}
0 & 1  \\
0 & 0  \\
\end{array} \right) \quad , \quad
Y=\left( \begin{array}{cc}
1 & 0  \\
0 & 0  \\
\end{array} \right) ,\nn \]
with the relation $[Y,X]=X$ as in Definition \eqref{definofCMH1}. On the other hand,  in the Hopf algebra $F(G_2)$ defined by
relations \eqref{fg2def}, the generators $a_n$ can be considered as the following functions on $G_2$:
\[ a_n(\fy)=\frac{1}{n!}\fy^{(n)}(0), \quad \forall \fy \in G_2  . \nn \]

The actions defined in Lemma \eqref{fg2actandg1coact}  can be realized in the following way. The exponentials $e^{tX}$ and $e^{tY}$,
as elements of the affine group $G_1$, are given by
\[e^{tX}=(1,t), \quad e^{tY}=(e^t,0). \nn \]
 The action $X\rhd a_n$ in lemma \eqref{fg2actandg1coact} can be identified as:
\begin{eqnarray*}
&(X\rhd a_n)(\fy)&=\frac{d}{dt}|_{{}_{t=0}}a_n(\fy \lhd e^{tX})\\
&&=\frac{d}{dt}|_{{}_{t=0}}a_n(\frac{\fy(x +t) -\fy(t)}{\dot{\fy}(t)})\\
&&=\frac{d}{dt}|_{{}_{t=0}}(\frac{\fy^{(n)}(t)}{n!\dot{\fy}(t)})\\
&&=\frac{d}{dt}|_{{}_{t=0}}(\frac{\fy^{(n)}(0)+t\fy^{(n+1)}(0)}{n!(\dot{\fy}(0) +t\ddot{\fy}(0))})\\
&&=\frac{d}{dt}|_{{}_{t=0}}(1/n!)(\fy^{(n)}(0)-t\ddot{\fy}(0)\fy^{(n)}(0) +t\fy^{(n+1)}(0)) \\
&&=(1/n!)(-\ddot{\fy}(0)\fy^{(n)}(0) +\fy^{(n+1)}(0)) \\
&&=[(n+1)a_{n+1}-2a_2a_n] (\fy),
\end{eqnarray*}
for any $\fy$ in $G_2$. This implies
\[X \rhd a_n= (n+1)a_{n+1}-2a_2a_n , \nn \]
as in Lemma  \eqref{fg2actandg1coact}. A similar computation will give the action $Y \rhd a_n$ as in Lemma  \eqref{fg2actandg1coact}.
 The coactions defined in Lemma  \eqref{fg2actandg1coact} can also be realized using the factorization
$G=G_1G_2$ \cite{maj3}.

Note that the above realization of actions and coactions is not necessary for the proof of Lemma \eqref{fg2actandg1coact}. It rather gives a good intuition
about where those formulas for actions and coactions come from.

\section{The supergroup $G^s=Diff^+(\mbr^{1,1})$ and its factorization} \label{superdiffgrp}

In this section, by replacing $\mbr$ by the supermanifold $\mbr^{1,1}$, we define a super version of the
group $G=Diff^+(\mbr)$, namely the supergroup $G^s=\dff$. Analogous to the non-super case,
we consider two sub supergroups $G^s_1$ and $G^s_2$
of $G^s$,  where  $G^s_1$ is the affine part of $G^s$. We establish the factorization $G^s=G^s_1G^s_2$.
This factorization will result in  a left action of $G^s_2$ on $G^s_1$ and a right action of $G^s_1$ on $G^s_2$.

For the general theory of supermanifolds we refer to \cite{ak6, r, v}. The supermanifold $\mbr^{1,1}$, the superline, is a super ringed space
 $S=( \mbr^1, \, \mathfrak{O}_S)$, where $ \mathfrak{O}_S$ is a sheaf of supercommutative $\mbr$-algebras over $\mbr$ defined,
 for each $U \subset \mbr$  open, by
\[\mathfrak{O}_S(U)= C^ \infty (U) \ot \Lambda^\bullet(\mbr) = C^ \infty (U) [\tta],  \nn\]
where $\tta$, the generator of the exterior algebra $\Lambda^\bullet(\mbr)$, is called the odd generator. We also normally denote the even
indeterminate  by $x$.
 Supermanifolds form a category where a morphism $f:S_1 \ra S_2$ is a morphism of the underlying super ringed spaces
\cite{ak6, r, v}.

A super Lie group is a group object in the category of supermanifolds. Alternatively, a super Lie group can be defined as a representable functor from the category of supermanifolds to the category of groups.
A typical example is $GL(p,q)= GL(\mbr^{p,q})$, the super general linear group of automorphisms of $\mbr^{p,q}$. This supergroup
has a matrix representation as follows. It is formed by matrices
\[
\left( \begin{array}{c|c}
A & B  \\
\hline
C & D  \\
\end{array} \right), \nn\]
where $A,D$ are, respectively, $p \times p$ and $q \times q$ invertible matrices consisting of
 even elements, and $C,D$ are, respectively,  $p \times q$ and $q \times p$ matrices consisting of
 odd elements.
As a special case, $GL(2,1)$ is formed by matrices
\[ \left( \begin{array}{cc|c}
a & b & x \\
c & d & y \\
\hline
z & w & e
\end{array} \right), \nn\]
where $a,b,c,d,e$ are even elements, $ad-bc \neq 0$, $e\neq 0$, and $x,y,z,w$ are odd elements.

\bde \label{supdiff}
The supergroup of orientation preserving diffeomorphisms of the super real line, $\mbr^{1,1}$, is defined as follows:
\[G^s=\dff=\{\Fy(x,\tta)=(A(x)+B(x)\tta \,,\, C(x)+D(x)\tta)\} , \nn\]
such that $A(x),D(x)$ are even, $ A(x),D(x) \in Diff(\mbr)$, $\dot{A}(x)>0,\, \dot{D}(x)>0$, and $B(x),C(x)$ are odd.
\ede

\bde \label{supaff}
The affine part of $G^s$, denoted by $G^s_1$ and also denoted by $\aff$, is defined by:
\[G^s_1 =\{\si(x, \tta)=(ax+b\tta +e \, ,\, cx+d\tta +f) \in G^s \}, \nn\]
such that $a, e, d$ are even, $ a , d > 0$,  and  $b, c, f$ are  odd.
An element $\si(x, \tta)=(ax+b\tta +e \, ,\, cx+d\tta +f)$ of $G^s_1$ can also be represented in the following way:
\[\si(x, \tta) =
\left( \begin{array}{c|c}
a & b \\
\hline
c & d
\end{array} \right)
\left( \begin{array}{c}
x   \\
\tta
\end{array} \right)
+
\left( \begin{array}{c}
e   \\
f
\end{array} \right)
\nn .\]
Therefore there exists the following representation of $G^s_1$ in $GL(2,1)$, which will be very useful for our purpose:
\[G^s_1 = \left\{ \si=
\left( \begin{array}{cc|c}
a & e & b \\
0 & 1 & 0 \\
\hline
c & f & d
\end{array} \right) \,| \, a, e, d~are~ even,~ b, c, f~are~ odd,~and~ a , d > 0\right\}. \nn\]
\ede

\bde \label{supnonaff}
A second super subgroup $G^s_2$ of $G^s$ is defined by:
\[G^s_2=\{\fy=(A(x)+B(x)\tta, C(x)+D(x)\tta) \in G^s \, |\, \fy(0,0)=0,\, J\fy(0,0)=1  \}. \nn\]
In other words
\beq\label{defofgs2}
G^s_2=\{\fy=(A(x)+B(x)\tta, C(x)+D(x)\tta) \in G^s \},
\eeq
where $ A(0)=B(0)=C(0)=\dot{C}(0)=0$, and $\dot{A}(0)=D(0)=1$.
\ede

In order to proceed to the factorization $G^s=G^s_1G^s_2$, for any $\Fy \in \dff$ we define
\beq \label{pi1}
\pi_1(\Fy)=(J\Fy(0,0),\Fy(0,0))  =
\left( \begin{array}{c|c}
\dot{A}(0) & B(0) \\
\hline
\dot{C}(0) & D(0)
\end{array} \right)
\left( \begin{array}{c}
x   \\
\tta
\end{array} \right)
+
\left( \begin{array}{c}
A(0)   \\
C(0)
\end{array} \right)
\in G^s_1
 ,
\eeq
and
\beq \label{pi2}
\begin{split}
\pi_2(\Fy)&= (J\Fy(0,0))^{-1} (\Fy(x,\tta)-\Fy(0,0))\\
&=
\left( \begin{array}{c|c}
\dot{A}(0) & B(0) \\
\hline
\dot{C}(0) & D(0)
\end{array} \right)^{-1}
(\Fy(x,\tta)-\Fy(0,0)) \in G^s_2.
\end{split}
\eeq
Here we have used the definition

\begin{displaymath}\nn
J\Fy(x,\tta):=
\left( \begin{array}{c|c}
\frac{\rnd\Fy_1}{ \rnd x} & - \frac{\rnd\Fy_1}{ \rnd \tta} \\
\hline
\frac{\rnd\Fy_2}{ \rnd x} & \frac{\rnd\Fy_2}{ \rnd \tta}
\end{array} \right)
=
\left( \begin{array}{c|c}
\dot{A}(x)+\dot{B}(x)\tta & B(x) \\
\hline
\dot{C}(x)+\dot{D}(x)\tta & D(x)
\end{array} \right),
\end{displaymath}
where $\Fy_1=A(x)+B(x)\tta$ and $\Fy_2=C(x)+D(x)\tta$ are, respectively, the even and odd components of $\Fy \in \dff$.
The operator $\frac{\rnd}{ \rnd x}$ is even and the operator $\frac{\rnd}{ \rnd \tta}$ is odd.
Also, the formula for the inverse supermatrix is
\[
\left( \begin{array}{c|c}
a & b \\
\hline
c & d
\end{array} \right)^{-1}
=\frac{1}{da}
\left( \begin{array}{c|c}
d+\frac{bc}{a} & -b \\
\hline
-c & a+\frac{cb}{d}
\end{array} \right)
. \nn \]

Now if we let $\si(x,\tta)=\pi_1(\Fy)$ and $\fy(x,\tta)=\pi_2(\Fy)$, it is clear that  for any $\Fy \in \dff$ 
\[\Fy=\si \fy, \nn \]
which proves the factorization
\[G^s=G^s_1 G^s_2 . \nn \]
Therefore, we have the following two natural actions.
The left action of $G^s_2$ on $G^S_1$, $G^s_2 \times G^s_1 \ra G^s_1$, defined by $\fy \rhd \si =\pi_1(\fy \si)$,
and the right action of $G^s_1$ on $G^s_2$, $G^s_2 \times G^s_1 \ra G^s_2$, defined by $\fy \lhd \si =\pi_2(\fy \si)$.

\section {The super Hopf algebras $U(\mfg^s_1)$ and $F(G^s_2)$}

In this section we assign to supergroups $G^s_1$ and $G^s_2$, the super Hopf algebras $U(\mfg^s_1)$ and $F(G^s_2)$, respectively. The super Hopf algebra $U(\mfg^s_1)$ is just
the universal enveloping algebra of the super Lie algebra $\mfg^s_1 =Lie(G^s_1)$.
As for $F(G^s_2)$, we define the coordinate functions
$a_n,b_n,c_n$ and $d_n$ on $G^s_2$. Then  $F(G^s_2)$ would be the corresponding
 {\it{Fa\`a di Bruno}} or rather {\it{comeasuring}} super Hopf algebra generated,
as a supercommutative superalgebra, by
  $a_n,b_n,c_n$ and $d_n$,  for which we will define the Hopf algebra structure.

\subsection{The super Hopf algebra $U(\mfg^s_1)$}

The super Lie algebra $\mfg^s_1 =Lie(G^s_1)$ is generated by three even generators:
\begin{eqnarray}
X= \left( \begin{array}{cc|c} 0 & 1 & 0 \\ 0 & 0 & 0 \\ \hline 0 & 0 & 0 \end{array} \right),   \quad
Y= \left( \begin{array}{cc|c} 1 & 0 & 0 \\ 0 & 0 & 0 \\ \hline 0 & 0 & 0 \end{array} \right), \quad
Z= \left( \begin{array}{cc|c} 0 & 0 & 0 \\ 0 & 0 & 0 \\ \hline 0 & 0 & 1 \end{array} \right) ,\nn
\end{eqnarray}
and three odd generators

\begin{eqnarray}
&U&= \left( \begin{array}{cc|c} 0 & 0 & 1 \\ 0 & 0 & 0 \\ \hline 0 & 0 & 0 \end{array} \right), \quad
V= \left( \begin{array}{cc|c} 0 & 0 & 0 \\ 0 & 0 & 0 \\ \hline -1 & 0 & 0 \end{array} \right),  \quad
W= \left( \begin{array}{cc|c} 0 & 0 & 0 \\ 0 & 0 & 0 \\ \hline 0 & -1 & 0 \end{array} \right).\nn
\end{eqnarray}

It is easy to check that the following \emph{super bracket} relations hold:
\[ [Y,X]=X, \quad [X,Z]=0, \quad [Y,Z]=0 , \nn \]
\[ [X,U]= 0, \quad [X,V]=-W, \quad [X,W]=0 , \nn \]
\[ [Y,U]= U, \quad [Y,V]= -V, \quad [Y,W]= 0 , \nn \]
\[ [Z,U]= -U, \quad [Z,V]= V, \quad [Z,W]= W , \nn \]
\[ [U,V]=-(Y+Z), \quad [U,W]= -X, \quad [V,W]= 0 . \nn \]
\beq\label{bracketres}
[X,X]=[Y,Y]=[Z,Z]=[U,U]= [V,V]= [W,W]= 0 .
\eeq

The super Hopf algebra $U(\mfg^s_1)$ is the universal enveloping algebra of $\mfg^s_1$.

\subsection {The super Hopf algebra $F(G^s_2)$}
The super Hopf algebra $F(G^s_2)$, as a supercommutative superalgebra,
is the super polynomial algebra $\mbr [a_n,b_n,c_n,d_n]$, generated by
two sets of even generators  $a_n,d_n, \, n \geq 0 $, $a_0=0, a_1=d_0=1$,
 and two sets of odd generators $b_n,c_n, \, n \geq 0 $, $b_0=c_0=c_1=0$,
where for any $\fy(x,\tta)=(A(x)+B(x)\tta \, , \, C(x)+D(x)\tta)$ in $G^s_2$ we have:
\[a_n(\fy)=(1/n!)A^{(n)}(0), \nn \]
\[b_n(\fy)=(1/n!)B^{(n)}(0), \nn \]
\[c_n(\fy)=(1/n!)C^{(n)}(0), \nn \]
\beq\label{anfy}
d_n(\fy)=(1/n!)D^{(n)}(0).
\eeq

To define the coproduct on  $F(G^s_2)$, analogous to the non-super case, we use the following formalism:
\beq\label{formlfrcomultinfg2s}
m \Del(a_n) (\fy \ot \fy') = a_n (\fy \circ \fy'), \quad \fy, \fy' \in G^s_2,
\eeq
and similar formulas for $b_n$, $c_n$ and $d_n$.
 Note that using this formalism provides us with the coassociativity property of $\Del$.

We need to study the composition of two elements of $G^s_2$.
Let $\fy(x,\tta)=(A(x)+B(x)\tta \, , \, C(x)+D(x)\tta)$ and $\fy'(x,\tta)=(A'(x)+B'(x)\tta \, , \, C'(x)+D'(x)\tta)$
be two elements of $G^s_2$. Let $f=A'(x)+B'(x)\tta$ and $g=C'(x)+D'(x)\tta$.
It is easy to check that the composition is given by
\[\fy \circ \fy'(x,\tta)=(A(f)+B(f)g \, , \, C(f)+D(f)g), \label{fyfy} \nn \]
where
\begin{eqnarray} \label{afbfg} \nn
&A(f)+B(f)g&=[A(A'(x)+ B(A'(x))C'(x)] + \\
&&[\dot{A}(A'(x))B'(x)+ B(A'(x))D'(x)-\dot{B}(A'(x))B'(x)C'(x)]\tta,\nn
\end{eqnarray}
and
\begin{eqnarray} \label{cfdfg} \nn
&C(f)+D(f)g&=[C(A'(x))+ D(A'(x))C'(x)] + \\
&&[\dot{C}(A'(x))B'(x)+ D(A'(x))D'(x)-\dot{D}(A'(x))B'(x)C'(x)]\tta. \nn
\end{eqnarray}

Now we proceed to comultiplications, starting with $\Del(a_n)$. It is not hard to show that
\[a_1(\fy \circ \fy')=1 = m(1\ot 1)(\fy \ot \fy'), \nn \]
\[a_2(\fy \circ \fy')=a'_2+a_2 = m(1\ot a_2 +a_2 \ot 1)(\fy \ot \fy'), \nn \]
\[a_3(\fy \circ \fy')=a'_3+2a_2a'_2+a_3 +b_1c'_2 = m(1\ot a_3+ 2a_2\ot a_2+a_3 \ot 1 +b_1\ot c_2)(\fy \ot \fy'). \nn \]
Therefore by \eqref{formlfrcomultinfg2s} we have:
\[\Del(a_1)= \Del(1)=1\ot 1 , \nn \]
\[\Del(a_2)=1\ot a_2 +a_2 \ot 1 , \nn \]
\[\Del(a_3)=1\ot a_3+ a_3 \ot 1+ 2a_2\ot a_2 +b_1\ot c_2, \nn \]
and more generally for $n\geq 1$,
\begin{multline*}
\Del(a_n)= \sum_{k=1}^n a_k \ot \sum_{l_1 + l_2 + \cdots + l_k =n} a_{l_1}a_{l_2} \cdots a_{l_k}+ \\
 \sum_{i=1}^{n} \sum_{k=1}^i b_k \ot (\sum_{l_1 + l_2 + \cdots + l_k =i} a_{l_1}a_{l_2} \cdots a_{l_k})c_{n-i} .
\end{multline*}

With a  similar method we have, for $\Del(b_n)$,
\[\Del(b_1)=1\ot b_1 +b_1 \ot 1 , \nn \]
\[\Del(b_2)=1\ot b_2 +b_2 \ot 1 +2a_2 \ot b_1 + b_1 \ot d_1 +b_1 \ot a_2 , \nn \]
\begin{multline*}
\Del(b_3)= 1\ot b_3+ b_3 \ot 1+ 2a_2\ot b_2 + 2a_2\ot a_2b_1 + 3a_3\ot b_1 + b_1\ot d_2   \\
 + b_1\ot a_2d_1 + b_2\ot d_1 + b_1\ot a_3 + 2b_2\ot a_2 - b_1\ot b_1c_2 ,
\end{multline*}
and in general
\begin{multline*}
\Del(b_n)= 1 \ot b_n + \sum_{i=1}^{n}  \sum_{k=1}^i (k+1)a_{k+1} \ot (\sum_{l_1 + l_2 + \cdots + l_k =i} a_{l_1}a_{l_2} \cdots a_{l_k}) b_{n-i} \\
+  \sum_{i=1}^{n} \sum_{k=1}^i b_k \ot ( \sum_{l_1 + l_2 + \cdots + l_k =i} a_{l_1}a_{l_2} \cdots a_{l_k})d_{n-i} \\
- \sum_{i=1}^{n} \left(b_1 \ot b_ic_{n-i} + \sum_{j=1}^{i}  \sum_{k=1}^j (k+1)b_{k+1} \ot (\sum_{l_1 + l_2 + \cdots + l_k =j} a_{l_1}a_{l_2}\cdots a_{l_k})b_{i-j}c_{n-i}\right).
\end{multline*}

For $\Del(c_n)$, $n\geq 1$, we have:
\begin{multline*}
\Del(c_n)= 1 \ot c_n + \sum_{k=1}^n c_k \ot \sum_{l_1 + l_2 + \cdots + l_k =n} a_{l_1}a_{l_2} \cdots a_{l_k}+ \\
\sum_{i=1}^{n} \sum_{k=1}^i d_k \ot (\sum_{l_1 + l_2 + \cdots + l_k =i} a_{l_1}a_{l_2} \cdots a_{l_k})c_{n-i} .
\end{multline*}
In particular for $n=2,3$:
\[\Del(c_2)=1\ot c_2 +c_2 \ot 1  , \nn \]
\[\Del(c_3)=1\ot c_3 +c_3 \ot 1 + 2c_2 \ot a_2 + d_1 \ot c_2. \nn \]

For  $\Del(d_n)$, $n\geq 1$, we have:
\begin{multline*}
\Del(d_n)= 1 \ot d_n + \sum_{i=1}^{n}  \sum_{k=1}^i (k+1)c_{k+1} \ot (\sum_{l_1 + l_2 + \cdots + l_k =i} a_{l_1}a_{l_2} \cdots a_{l_k}) b_{n-i} \\
+  \sum_{i=1}^{n} \sum_{k=1}^i d_k \ot ( \sum_{l_1 + l_2 + \cdots + l_k =i} a_{l_1}a_{l_2} \cdots a_{l_k})d_{n-i}\\
- \sum_{i=1}^{n} \left(d_1 \ot b_ic_{n-i} + \sum_{j=1}^{i}  \sum_{k=1}^j (k+1)d_{k+1} \ot (\sum_{l_1 + l_2 + \cdots + l_k =j} a_{l_1}a_{l_2}\cdots a_{l_k})b_{i-j}c_{n-i}\right).
\end{multline*}
In particular for $n=1,2,3$:
\[\Del(d_1)=1\ot d_1 +d_1 \ot 1 , \nn \]
\[\Del(d_2)=1\ot d_2 +d_2 \ot 1 +2c_2 \ot b_1 + d_1 \ot d_1 +d_1 \ot a_2 , \nn \]
\begin{multline*}
\Del(d_3)= 1\ot d_3+ d_3 \ot 1+ 2c_2\ot b_2 + 2c2\ot a_2b_1 + 3c_3\ot b_1 + d_1\ot d_2  \\
 + d_1\ot a_2d_1 + d_2\ot d_1 + d_1\ot a_3 + 2d_2\ot a_2 - d_1\ot b_1c_2 .
\end{multline*}

We extend $\Del$ linearly to $F(G_2)$ via the relation
\[ \Del (ab):= \Del (a)\Del (b), \forall a,b \in F(G_2)  ,\]
where multiplication on the right hand side operates in the super sense. We also set $\epn(1)=1$, and define $\epn$ to be equal to zero on all other generators. This defines a super bialgebra structure on $F(G_2)$.

To prove that $F(G^s_2)$ is a super Hopf algebra, the only missing data is the antipode map. Analogous to the non-super case, it suffices to define an anti-algebra and coalgebra map $S:F(G^s_2) \to F(G^s_2)$ satisfying
\beq \nn
S(a) (\fy) = a(\fy^{-1}), \quad \fy \in G^s_2, \quad a= a_n,b_n,c_n,d_n .
\eeq
We shall do this  in an inductive fashion at the end of the next section, where we define the actions and coactions between $F(G^s_2)$ and $U(\mfg^s_1)$. The reason behind this is that we want to use the fact that the antipode, if it exists, should interact, {\it{in a nice way}}, with those actions and coactions, in the sense of relation \eqref{bmform} below \cite{bh10}.

\section{Actions and coactions} \label{actcoactbicrsprdct}
In this section we prove that $F(G^s_2)$ is a super left $U(\mfg^s_1)$-module algebra
and $U(\mfg^s_1)$ is a super right  $F(G^s_2)$-comodule coalgebra.

For $Y= \left( \begin{array}{cc|c} 1 & 0 & 0 \\ 0 & 0 & 0 \\ \hline 0 & 0 & 0 \end{array} \right) \in \mfg^s_1$ we have
\[e^{tY}= \left( \begin{array}{cc|c} e^t & 0 & 0 \\ 0 & 1 & 0 \\ \hline 0 & 0 & 1 \end{array} \right)=(e^t x, \tta)  \in G^s_1,\quad (t~even) .
 \nn\]
The action $Y\rhd a_n$ can be realized by
\[(Y\rhd a_n)(\fy)=\frac{d}{dt}|_{{}_{t=0}}a_n(\fy \lhd e^{tY}),  \nn\]
for any $\fy=(A(x)+B(x)\tta, C(x)+D(x)\tta) \in G^s_2$.
Thus, by computing $\fy \lhd e^{tY}=\pi_2(\fy e^{tY})$,
we obtain
\[a_n(\fy \lhd e^{tY})=(1/n!) e^{(n-1)t}A^{(n)}(0), \nn\]
and
\[(Y\rhd a_n)(\fy)=\frac{d}{dt}|_{{}_{t=0}}a_n(\fy \lhd e^{tY})= (n-1)a_n(\fy). \nn \]
This implies that
\[Y\rhd a_n=(n-1)a_n. \nn \]
In the same way we have
\[Y\rhd b_n=(n-1)b_n ,\nn \]
\[Y\rhd c_n=nc_n, \nn \]
\[Y\rhd d_n=nd_n. \nn \]

Using the same method we can derive all other actions. In fact we have the following lemma:
\ble  \label{FG2isUg1modalg}
Let us define the actions of $X, Y, Z, U, V, W$ on $a_n, b_n, c_n, d_n$ by the following relations:
\[X\rhd a_n=(n+1)a_{n+1}-2a_na_2 -b_1c_n, \quad X\rhd b_n=(n+1)b_{n+1}-2b_na_2 -b_1d_n ,\nn \]
\[X\rhd c_n=-2c_2 a_n+(n+1)c_{n+1}-c_nd_1  , \quad X\rhd d_n=-2c_2 b_n+(n+1)d_{n+1}-d_nd_1, \nn \]
\[Y\rhd a_n=(n-1)a_n, \quad Y\rhd b_n=(n-1)b_n ,\quad Y\rhd c_n=nc_n, \quad Y\rhd d_n=nd_n, \nn \]
\[Z\rhd a_n=0, \quad Z\rhd b_n=b_n ,\quad Z\rhd c_n=-c_n, \quad Z\rhd d_n=0, \nn \]
\[U\rhd a_n=c_n, \quad U\rhd b_n=-(n+1)a_{n+1}+d_n ,\quad U\rhd c_n=0, \quad U\rhd d_n=(n+1)c_{n+1}, \nn \]
\[V\rhd a_n=b_{n-1}, \quad V\rhd b_n=0 ,\quad V\rhd c_n=a_n-d_{n-1}, \quad V\rhd d_n=+b_n, \nn \]
\[W\rhd a_n=-b_1a_n+b_n, \quad W\rhd b_n=-b_1b_n , \nn\]
\beq\label{actionofsxyzuvw}
 W\rhd c_n=d_1a_n-d_{n}, \quad W\rhd d_n=d_1b_n,
\eeq
and extend those actions to an action of  $U(\mfg^s_1)$  on  $F(G^s_2)$ such that
\beq\label{ghrhda}
(gh) \rhd a = g\rhd (h\rhd a),
\eeq
\beq\label{hrhdab}
h \rhd (ab):=(-1)^{|a||h_{(2)}|} \quad (h_{(1)}\rhd a)(h_{(2)}\rhd b),
\eeq
for $a, b$ in $F(G^s_2)$ and  $g, h$ in  $U(\mfg^s_1)$. Then $F(G^s_2)$ is a super left $U(\mfg^s_1)$-module algebra.
\ele

\bpf
It is enough to show that this action is consistent with the  bracket relations \eqref{bracketres}.
We check this just for some of the bracket relations. The rest would be the same:\\
\beq \nn
\begin{split}
(YX)\rhd a_n &\overset{\eqref{ghrhda}}{ =} Y\rhd (X\rhd a_n)   \overset{\eqref{actionofsxyzuvw}}{ =}  Y \rhd ((n+1)a_{n+1} -2 a_na_2 -b_1c_n )\\
& = (n+1)(Y \rhd a_{n+1}) -2  Y \rhd (a_na_2) -  Y \rhd (b_1c_n ) \\
&\overset{\eqref{hrhdab}}{=} n(n+1)a_{n+1} -2(a_n(Y\rhd a_2) +(Y\rhd a_n)a_2) - (b_1(Y\rhd c_n) +(Y\rhd b_1)c_n) \\
&= (n^2+n)a_{n+1} -2n a_n a_2 -nb_1 c_n,\\
(XY)\rhd a_n &= X\rhd (Y\rhd a_n)= X \rhd ((n-1)a_n) \\
&= (n-1)((n+1)a_{n+1} -2a_na_2 -b_1 c_n)\\
&=(n^2 - 1)a_{n+1} -2(n-1)a_na_2 -(n-1)b_1 c_n .
\end{split}
\eeq
Therefore,
\[[Y,X]\rhd a_n \overset{}{=} (YX -XY)\rhd a_n =(n+1)a_{n+1} -2a_na_2 -b_1c_n  \overset{\eqref{actionofsxyzuvw}}{ =} X\rhd a_n, \nn\]
which is consistent with the relation $[Y,X]= X$ of \eqref{bracketres}. Next we check a bracket involving odd generators:
\beq \nn
\begin{split}
(UW)\rhd d_n &= U\rhd (W\rhd d_n)= U \rhd (d_1b_n)\\
&= (d_1(U \rhd b_n) + (U\rhd d_1)b_n)= d_1 (-(n+1)a_{n+1}+d_n) + 2c_2 b_n \\
&=  -(n+1) d_1a_{n+1}+  d_1d_n + 2c_2 b_n,\\
(WU)\rhd d_n &= W \rhd (U\rhd d_n)= W \rhd ((n+1)c_{n+1}) = (n+1)d_1a_{n+1}- (n+1)d_{n+1}  .
\end{split}
\eeq
Therefore,
\[ [U,W]\rhd d_n =(UW+WU) \rhd d_n =   2c_2 b_n - (n+1)d_{n+1} +  d_1d_n =(-X) \rhd d_n, \nn\]
which agrees with the relation $[U,W]= -X$.\\
 \epf

By using almost the same method we can find the coactions and prove the following lemma.
 Let us denote the right coaction of $F(G^s_2)$ on  $U(\mfg^s_1)$, by
$\nabla_r : U(\mfg^s_1) \ra U(\mfg^s_1) \ot F(G^s_2)$.

\ble \label{Ug1Fg2comodcoalg}
Let us define the coactions of  $F(G^s_2)$ on generators $X,Y,Z,U,V,W$ of  $U(\mfg^s_1)$ by
\[\nabla_r(X)=  2Y\ot a_2 +  X\ot 1 + Z \ot d_1 + U\ot b_1 + 2V\ot c_2 , \quad \nabla_r(Y)= Y\ot 1 , \nn \]
\[\nabla_r(Z)= Z\ot 1 , \quad \nabla_r(U)= U\ot 1 , \quad \nabla_r(V)= V\ot 1, \nn \]
\beq\label{coactionsonxyzuvw}
\nabla_r(W)= Y\ot b_1 +V\ot d_1 + W\ot 1 ,
\eeq
and extend them to  $U(\mfg^s_1)$ such that
\beq\label{coactofgh}
\nabla_r(gh)= (-1)^{|h^{(1)}|(|{g_{(1)}}^{(2)}|+|g_{(2)}|)} \quad  {g_{(1)}}^{(1)} h^{(1)} \ot {g_{(1)}}^{(2)}(g_{(2)} \rhd h^{(2)}) ,
\eeq
for all $g$ and $h$ in $U(\mfg^s_1)$. Then $U(\mfg^s_1)$ is a super right  $F(G^s_2)$-comodule coalgebra.
\ele

\bpf
It is straightforward to check the coaction property, $(id\ot \Del)\nabla_r (h) = (\nabla_r \ot id)\nabla_r (h)$, for all $h$ in  $U(\mfg^s_1)$, in other words:
\[{h^{(1)}}\ot  {h^{(2)}}_{(1)} \ot  {h^{(2)}}_{(2)} = {h^{(1)(1)}}\ot {h^{(1)(2)}}\ot  {h^{(2)}}. \nn\]
We prove that this coaction is consistent with the bracket
relations \eqref{bracketres}.  We verify this only for one of the purely odd cases. The rest are similar. By formula \eqref{coactofgh} we have

\beq \nn
\begin{split}
\nabla_r(VW)&= (-1)^{|W^{(1)}|(|{V_{(1)}}^{(2)}|+|V_{(2)}|)} \quad  {V_{(1)}}^{(1)} W^{(1)} \ot {V_{(1)}}^{(2)}(V_{(2)} \rhd W^{(2)}) \\
&= [(-1)^{|W^{(1)}||V|)} \quad   W^{(1)} \ot (V \rhd W^{(2)})] \\
   &\quad  +  [ (-1)^{|W^{(1)}||V^{(2)}|} \quad  V^{(1)} W^{(1)} \ot V^{(2)} W^{(2)} ] \\
&= [(-1)^{|W^{(1)}|} \quad   W^{(1)} \ot (V \rhd W^{(2)})]+  [   V^{(1)} W^{(1)} \ot W^{(2)} ]  \\
&= [   Y\ot( V \rhd b_1) - V\ot (V \rhd  d_1) - W\ot  (V \rhd 1)   ] \\
   &\quad +  [  VY \ot b_1 + V^2 \ot d_1 + VW \ot 1]  \\
&= [  - V\ot b_1   ] +  [  VY \ot b_1 + V^2 \ot d_1 + VW \ot 1  ,
\end{split}
\eeq
and
\beq \nn
\begin{split}
\nabla_r(WV)&= (-1)^{|V^{(1)}|(|{W_{(1)}}^{(2)}|+|W_{(2)}|)} \quad  {W_{(1)}}^{(1)} V^{(1)} \ot {W_{(1)}}^{(2)}(W_{(2)} \rhd V^{(2)})\\
&= [(-1)^{|V^{(1)}||W|)} \quad   V^{(1)} \ot (W \rhd V^{(2)})] \\
   &\quad +  [ (-1)^{|V^{(1)}||W^{(2)}|} \quad  W^{(1)} V^{(1)} \ot W^{(2)} V^{(2)} ]\\
&= [-   V \ot (0)] +  [ (-1)^{|W^{(2)}|} \quad  W^{(1)} V \ot W^{(2)} ]  \\
&= -YV\ot b_1 + V^2 \ot d_1 + WV \ot 1 .\\
\end{split}
\eeq
Thus,
\beq \nn
\begin{split}
\nabla_r([V,W])& = \nabla_r(VW +WV)= \nabla_r(VW)+ \nabla_r(WV) \\
&= [V,Y] \ot b_1 + [V,V]\ot d_1 + [V,W] \ot 1 - v \ot b_1 = 0,
\end{split}
\eeq
which agrees with $[V,W] = 0$ of relations \eqref{bracketres}.

We also leave it to the reader to check that $U(\mfg^s_1)$ is a right  $F(G^s_2)$-comodule coalgebra, i.e.,  for all $h$ in $U(\mfg^s_1)$,
\[{h^{(1)}}_{(1)}\ot {h^{(1)}}_{(2)}\ot  h^{(2)}= (-1)^{|{{h_{(2)}}^{(1)}}||{{h_{(1)}}^{(2)}}|} {{h_{(1)}}^{(1)}}\ot {{h_{(2)}}^{(1)}}\ot {{h_{(1)}}^{(2)}}{{h_{(2)}}^{(2)}} . \nn \]

 \epf

\subsection*{Antipode for $F(G^s_2)$}
Now we prove that $F(G^s_2)$ is a super Hopf algebra, by defining an anti-algebra and coalgebra map $S:F(G^s_2) \to F(G^s_2)$ satisfying the antipode property or the following identity:
\beq \label{formlfrantipode}
S(a) (\fy) = a(\fy^{-1}), \quad \fy \in G^s_2, \quad a= a_n,b_n,c_n,d_n .
\eeq
We do this, inductively, by defining $S$ on generators  $a_n,b_n,c_n,d_n$, and then extend it linearly to $F(G^s_2)$ via the following relations:
\[S(ab) := (-1)^{|a||b|} S(b)S(a) , \]
\[\Del (S(a)) := (-1)^{|a_{(1)}||a_{(2)}|} S(a_{(2)}) \ot S(a_{(1)}) . \]

Let us define, first, $S(1)=1$, and $S(a)= -a$, for $a= a_2,b_1,c_2,d_1$. It is obvious that the antipode property holds for these elements. Now suppose relation \eqref{formlfrantipode} is true for all elements $a_i,b_i,c_i,d_i$, $i \leq n$.  The actions of $X$ defined in relations \eqref{actionofsxyzuvw} in Lemma \ref{FG2isUg1modalg} allow us to  write the higher degree elements, $a_{n+1},b_{n+1},c_{n+1},d_{n+1}$, in terms of some other elements of lower degree, $a_i,b_i,c_i,d_i$, $i \leq n$. Therefore, to finish the process it is enough to prove that
\beq \label{formlfrantipodeaction}
S(X \rhd a) (\fy) = (X \rhd a)(\fy^{-1}), \quad \fy \in G^s_2, \,\, a=a_i,b_i,c_i,d_i, \, i \leq n .
\eeq
To prove this identity we need the following three lemmas.

\ble \label{bmlemma} \cite{bh10}
The antipode $S:F(G^s_2) \to F(G^s_2)$, if it exists, should satisfy the following relation:
\beq \label{bmform}
S(g \rhd a) = (g^{(1)} \rhd S(a)) S(g^{(2)}), \quad g \in U(\mfg^s_1), a \in F(G^s_2)
\eeq
\ele

\ble \label{seclem}
For any $\fy= (A(x)+B(x)\tta ,\, C(x)+D(x)\tta)$ in $G^s_2$, one can prove
\beq \label{seclemform}
\frac{d}{dt}|_{{}_{t=0}} \left( e^{t X^{(1)}} X^{(2)}(\fy) \right) =  \frac{d}{dt}|_{{}_{t=0}} \left( \fy \rhd e^{t X} \right) .
\eeq
\ele

\bpf
From the coactions defined in Lemma \ref{Ug1Fg2comodcoalg} we have:
\begin{align*}
\text{LHS} &=\frac{d}{dt}|_{{}_{t=0}} \left( e^{t Y} 2a_2(\fy) + e^{t X} 1(\fy) + e^{t Z} d_1(\fy) - e^{t U} b_1(\fy) - e^{t V} 2c_2(\fy) \right)\\
&=\frac{d}{dt}|_{{}_{t=0}} ( (e^tX, \tta) 2a_2 + (X+t, \tta) + (X, e^t\tta) d_1 \\
& \qquad \qquad \qquad \qquad -  (X+t\tta, \tta)b_1 - (X, -tX+\tta) 2c_2 )\\
&= (X, 0) 2a_2 + (1, 0) + (0, \tta) d_1 -  (\tta, 0)b_1 - (0, -X) 2c_2 \\
&= (2a_2 X+1+b_1 \tta , 2c_2X+d_1 \tta )
\end{align*}
From the discussion in Section \ref{superdiffgrp} and formula \eqref{pi1}, we have
\begin{align*}
\text{LHS} &=\frac{d}{dt}|_{{}_{t=0}} \left( \pi_1(\fy(X+t, \tta))  \right)
=\frac{d}{dt}|_{{}_{t=0}} \left(
\left( \begin{array}{c|c}
\dot{A}(t) & B(t) \\
\hline
\dot{C}(t) & D(t)
\end{array} \right)
\left( \begin{array}{c}
x   \\
\tta
\end{array} \right)
+
\left( \begin{array}{c}
A(t)   \\
C(t)
\end{array} \right) \right) \\
&=\frac{d}{dt}|_{{}_{t=0}} \left(
\dot{A}(t)X +A(t)+B(t)\tta , \dot{C}(t)X +C(t)+D(t)\tta  \right) \\
&=\left(
\ddot{A}(0)X +\dot{A}(0)+ \dot{B}(0)\tta , \ddot{C}(0)X +\dot{C}(0)+\dot{D}(0)\tta  \right)\\
&= (2a_2 X+1+b_1 \tta , 2c_2X+d_1 \tta )\\
&= \text{RHS}
\end{align*}
\epf

\ble \label{formfact}
If $G=G_1G_2$ is a factorisation of the group $G$, then, for any $\vfy_1$, $\vfy_2$ in $G_2$ and $\si$ in $G_1$, one has 
\beq
(\vfy_1 \vfy_2) \lhd \si = (\vfy_1 \lhd (\vfy_2 \rhd \si)) \, (\vfy_2 \lhd \si)
\eeq
\ele

Now we prove the identity \eqref{formlfrantipodeaction}, $S(X \rhd a) (\fy) = (X \rhd a)(\fy^{-1})$.
\begin{proof}
\begin{align*}
S(X \rhd a) (\fy) & \overset{Lemma \, \ref{bmlemma}}{=} (X^{(1)} \rhd S(a))(\fy)\, S(X^{(2)})(\fy)\\
&\overset{}{=} \frac{d}{dt}|_{{}_{t=0}} \left( S(a) ( \fy \lhd e^{t X^{(1)}} )\right) \, S(X^{(2)})(\fy)\\
&\overset{Induction}{=} \frac{d}{dt}|_{{}_{t=0}} \left( a ( ( \fy \lhd e^{t X^{(1)}} )^{-1}) \right) \, X^{(2)}(\fy^{-1})\\
&\overset{Lemma \, \ref{formfact}}{=} \frac{d}{dt}|_{{}_{t=0}} \left( a  ( \fy^{-1} \lhd (\fy \rhd e^{t X^{(1)}} )) \right) \, X^{(2)}(\fy^{-1})\\
&\overset{Lemma \, \ref{seclem}}{=} \frac{d}{dt}|_{{}_{t=0}} \left( a  ( \fy^{-1} \lhd e^{t X} ) \right) \\
&\overset{}{=}(X \rhd a)(\fy^{-1})
\end{align*}

\end{proof}

\section{Compatibilities and the super Hopf algebra  $\mh ^s_1$}

The following proposition is the super analogue of  Lemma \eqref{cptbltycondnonsup}. It gives the compatibility conditions  to construct a
bicrossproduct super Hopf algebra.
To complete the construction of the  bicrossproduct super Hopf algebra $F(G^s_2) \lrbicross U(\mfg^s_1)$ one needs to check  these compatibility conditions between the actions  and coactions introduced in Lemmas
\eqref{FG2isUg1modalg} and \eqref{Ug1Fg2comodcoalg} in the last section.

\bpr \label{supercptbltycondnonsup}
Let $A$ and $H$ be two super Hopf algebras such that $A$ is a left $H$-module algebra, and $H$ is a right $A$-comodule coalgebra.
Let furthermore these structures satisfy the follow compatibility conditions:
\beq\label{delhrhda}
\Del(h \rhd a)= (-1)^{|a_{(1)}|(|{h_{(1)}}^{(2)}|+|h_{(2)}|)} \quad  {h_{(1)}}^{(1)} \rhd a_{(1)} \ot {h_{(1)}}^{(2)} (h_{(2)} \rhd a_{(2)}) ,
\eeq
\[\epn(h\rhd a)=\epn(h) \epn(a)    , \nn \]
\beq\label{nablgh}
\nabla_r(gh)= (-1)^{|h^{(1)}|(|{g_{(1)}}^{(2)}|+|g_{(2)}|)} \quad  {g_{(1)}}^{(1)} h^{(1)} \ot {g_{(1)}}^{(2)}(g_{(2)} \rhd h^{(2)}) ,
\eeq
\[ \nabla_r(1)=1\ot 1     , \nn \]
\beq\label{cpcondbicrsprdct}
(-1)^{|h_{(1)}||{h_{(2)}}^{(1)}|+|a||{h_{(2)}}^{(2)}| } \quad   {h_{(2)}}^{(1)} \ot ({h_{(1)}}\rhd a){h_{(2)}}^{(2)} ={h_{(1)}}^{(1)} \ot {h_{(1)}}^{(2)} (h_{(2)} \rhd a) ,
\eeq
 for any $a,b$ in $A$ and $g,h$ in $H$, where we have denoted the actions by, $h \rhd a$ and the coactions by $\nabla_r (h)= h_{(0)} \ot  h_{(1)}$.
Then the super vector space $A\ot H$ can be equipped with a super Hopf algebra structure as follows:
\beq\label{aothbotg}
(a\ot h)(b\ot g)= (-1)^{|h_{(2)}||b|} \quad  a(h_{(1)} \rhd b) \ot h_{(2)} g ,
\eeq
\beq\label{delaotb}
\Del(a \ot h)=(-1)^{|{h_{(1)}}^{(1)}||{a_{(2)}}| } \quad  a_{(1)} \ot {h_{(1)}}^{(1)} \ot a_{(2)}{h_{(1)}}^{(2)}\ot h_{(2)}  ,
\eeq
\[\epn(a \ot h)= \epn(a)\epn(h)     , \nn \]
\[S(a \ot h)= (-1)^{|h^{(1)}||a|} \quad (1\ot S(h^{(1)}))(S(a h^{(2)}) \ot 1)     ,\nn \]
for any $a,b$ in $A$ and $g,h$ in $H$.
 We call this super Hopf algebra the left-right bicrossproduct super Hopf algebra $A \lrbicross H$.
\epr

\bpf
We prove that $\Del$, defined in \eqref{delaotb}, is an algebra map.
First we have:
\beq\label{yyy}
\begin{split}
\Del &((a\ot h)\cdot (b\ot g)) \overset{\eqref{aothbotg}}{=} (-1)^{^{\alf_0}} \quad \Del(a(h_{(1)} \rhd b)\ot h_{(2)}g)\\
&\overset{\eqref{delaotb}}{=} (-1)^{^{\alf_0+\alf_1}}\quad (a(h_{(1)}\rhd b))_{(1)} \ot {( h_{(2)}g)_{(1)}}^{(1)} \ot (a(h_{(1)} \rhd b))_{(2)}{( h_{(2)}g)_{(1)}}^{(2)}    \\
   &\quad \quad \ot ( h_{(2)}g)_{(2)} \\
&\overset{}{=} (-1)^{^{\alf_0+\alf_1+\alf_2+\alf_3}}\quad a_{(1)}(h_{(1)}\rhd b)_{(1)} \ot {( h_{(2)(1)}g_{(1)})}^{(1)}  \\
   &\quad \quad \ot a_{(2)}(h_{(1)} \rhd b)_{(2)}{( h_{(2)(1)}g_{(1)})}^{(2)}\ot  h_{(2)(2)}g_{(2)}\\
&\overset{\eqref{delhrhda},\eqref{nablgh}}{=} (-1)^{^{\alf_0+\alf_1+\alf_2+\alf_3+\alf_4+\alf_5}} \quad a_{(1)}({h_{(1)(1)}}^{(1)}\rhd b_{(1)}) \ot {( {h_{(2)(1)(1)}}^{(1)}{g_{(1)}}^{(1)})}  \\
   &\quad \quad  \ot a_{(2)}{h_{(1)(1)}}^{(2)}(h_{(1)(2)}\rhd b_{(2)})( {h_{(2)(1)(1)}}^{(2)}(h_{(2)(1)(2)} \rhd {g_{(1)}}^{(2)})) \ot  h_{(2)(2)}g_{(2)}\\
  &\overset{coassociativity}{=} (-1)^{^{\alf_0+\alf_1+\alf_2+\alf_3+\alf_4+\alf_5}} \quad a_{(1)}({h_{(1)}}^{(1)}\rhd b_{(1)}) \ot {( {h_{(3)}}^{(1)}{g_{(1)}}^{(1)})}   \\
&\quad \quad \ot a_{(2)}{h_{(1)}}^{(2)}(h_{(2)}\rhd b_{(2)})( {h_{(3)}}^{(2)}(h_{(4)} \rhd {g_{(1)}}^{(2)})) \ot  h_{(5)}g_{(2)},
\end{split}
\eeq
where
\beq \nn
\begin{split}
\alf_0 &=|h_{(2)}||b| = (|{h_{(2)(1)(1)}}|+||{h_{(2)(1)(2)}}|+||{h_{(2)(2)}}|) (|b_{(1)}|+|b_{(2)}|)\\
&=  (|h_{(3)}|+||h_{(4)}|+||h_{(5)}|) (|b_{(1)}|+|b_{(2)}|) \\
\alf_1&= |{( h_{(2)}g)_{(1)}}^{(1)}| | (a(h_{(1)} \rhd b)_{(2)}|\\
&= |{( {h_{(2)(1)(1)}}^{(1)}{g_{(1)}}^{(1)})}| |a_{(2)}{h_{(1)(1)}}^{(2)}(h_{(1)(2)}\rhd b_{(2)})|\\
&= (|{{h_{(2)(1)(1)}}^{(1)}|+|{g_{(1)}}^{(1)}}|) (|a_{(2)}|+|{h_{(1)(1)}}^{(2)}|+|h_{(1)(2)}|+|b_{(2)}|)\\
&=  (|{{h_{(3)}}^{(1)}|+|{g_{(1)}}^{(1)}}|) (|a_{(2)}|+|{h_{(1)}}^{(2)}|+|h_{(2)}|+|b_{(2)}|) \\
&=   |{g_{(1)}}^{(1)}||a_{(2)}|+  |{g_{(1)}}^{(1)}||{h_{(1)}}^{(2)}|+  |{g_{(1)}}^{(1)}||h_{(2)}|+  |{g_{(1)}}^{(1)}||b_{(2)}| \\
& \quad  +    |{h_{(3)}}^{(1)}||a_{(2)}|+ |{h_{(3)}}^{(1)}||{h_{(1)}}^{(2)}|+ |{h_{(3)}}^{(1)}||h_{(2)}|+ |{h_{(3)}}^{(1)}||b_{(2)}| \\
\alf_2&= |a_{(2)}||(h_{(1)}\rhd b)_{(1)}|= |a_{(2)}| |({h_{(1)(1)}}^{(1)}\rhd b_{(1)})|\\
&= |a_{(2)}| (|{h_{(1)(1)}}^{(1)}|+|b_{(1)}|)\\
&=  |a_{(2)}| (|{h_{(1)}}^{(1)}|+|b_{(1)}|) =  |a_{(2)}||{h_{(1)}}^{(1)}|+ |a_{(2)}||b_{(1)}| \\
\alf_3&= | h_{(2)(2)}||g_{(1)}|  =   | h_{(5)}||g_{(1)}|  \\
\alf_4&= |b_{(1)}|(|(h_{(1)(2)}| +|{h_{(1)(1)}}^{(2)}|) \\
&= |b_{(1)}|(|(h_{(2)}| +|{h_{(1)}}^{(2)}|) = |b_{(1)}||(h_{(2)}| + |b_{(1)}||{h_{(1)}}^{(2)}| \\
\alf_5&= |{g_{(1)}}^{(1)}|( | h_{(2)(1)(2)}|+|{h_{(2)(1)(1)}}^{(2)}|)  \\
&= |{g_{(1)}}^{(1)}|( | h_{(4)}|+|{h_{(3)}}^{(2)}|)  =  |{g_{(1)}}^{(1)}| | h_{(4)}|+  |{g_{(1)}}^{(1)}||{h_{(3)}}^{(2)}|
\end{split}
\eeq

On the other hand we have
\beq \nn
\begin{split}
\Del &(a\ot h)\cdot \Del(b\ot g)\\
&\overset{\eqref{delaotb}}{=} (-1)^{^{\bt_0}} \quad (a_{(1)} \ot {h_{(1)}}^{(1)} \ot a_{(2)}{h_{(1)}}^{(2)}\ot h_{(2)}) \cdot \\
   &\quad \quad (b_{(1)} \ot {g_{(1)}}^{(1)} \ot b_{(2)}{g_{(1)}}^{(2)}\ot g_{(2)}) \\
&\overset{}{=} (-1)^{^{\bt_0 +\bt_1}} \quad  (a_{(1)} \ot {h_{(1)}}^{(1)}) \cdot \\
   &\quad \quad (b_{(1)} \ot {g_{(1)}}^{(1)})  \ot (a_{(2)}{h_{(1)}}^{(2)}\ot h_{(2)}) \cdot ( b_{(2)}{g_{(1)}}^{(2)}\ot g_{(2)}) \\
&\overset{\eqref{aothbotg}}{=} (-1)^{^{\bt_0 +\bt_1 +\bt_2}} \quad   a_{(1)}({{h_{(1)}}^{(1)}}_{(1)} \rhd b_{(1)}) \ot {{h_{(1)}}^{(1)}}_{(2)}  {g_{(1)}}^{(1)}  \\
   &\quad \quad \ot   a_{(2)}{h_{(1)}}^{(2)}(h_{(2)(1)} \rhd  b_{(2)}{g_{(1)}}^{(2)}) \ot h_{(2)(2)} g_{(2)}
\end{split}
\eeq
where
$ \bt_0= |{h_{(1)}}^{(1)}||{a_{(2)}}| + |{g_{(1)}}^{(1)}||{b_{(2)}}|$,
$ \bt_1= |(b_{(1)} \ot {g_{(1)}}^{(1)})| |(a_{(2)}{h_{(1)}}^{(2)}\ot h_{(2)})|$, and
$ \bt_2=|b_{(1)}| |{{h_{(1)}}^{(1)}}_{(2)}| + | b_{(2)}{g_{(1)}}^{(2)}||h_{(2)(2)}|$.

Continuing this computation, using the standard sign rules for superalgebras, coalgebras and Hopf algebras and coassociativity of $\Del_H$,
we obtain
\beq \label{xxx}
\begin{split}
\Del (a\ot h)\cdot \Del(b\ot g) &\overset{}{=}
(-1)^{^{\bt_0 +\bt_1 +\bt_2+\bt_3-\bt_4-\bt_5-\bt_6+\bt_7-\bt_8+\bt_9-\bt_{10}+\bt_{11}-\bt_{12}+\bt_{13}+\bt_{14}+\bt_{15}}} \\
  &\quad \quad (a_{(1)}{{h_{(1)}}^{(1)}} \rhd b_{(1)}\ot {{h_{(3)}}^{(1)}} {g_{(1)}}^{(1)}   \\
 &\quad \quad\ot  a_{(2)} {h_{(1)}}^{(2)} (h_{(2)} \rhd  b_{(2)}){h_{(3)}}^{(2)}(h_{(4)} \rhd {g_{(1)}}^{(2)}) \ot h_{(5)} g_{(2)})
\end{split}
\eeq
where at the end
\beq\label{zzz}  \nn
\begin{split}
\bt_0 &+\bt_1 +\bt_2+\bt_3-\bt_4-\bt_5-\bt_6+\bt_7-\bt_8+\bt_9-\bt_{10}+\bt_{11}-\bt_{12}+\bt_{13}+\bt_{14}+\bt_{15}\\
&= (|h_{(3)}|+||h_{(4)}|+||h_{(5)}|) (|b_{(1)}|+|b_{(2)}|)  \\
  &\quad + |{g_{(1)}}^{(1)}||a_{(2)}|+  |{g_{(1)}}^{(1)}||{h_{(1)}}^{(2)}|+  |{g_{(1)}}^{(1)}||h_{(2)}|+  |{g_{(1)}}^{(1)}||b_{(2)}| \\
  &\quad  +    |{h_{(3)}}^{(1)}||a_{(2)}|+ |{h_{(3)}}^{(1)}||{h_{(1)}}^{(2)}|+ |{h_{(3)}}^{(1)}||h_{(2)}|+ |{h_{(3)}}^{(1)}||b_{(2)}|  \\
  &\quad  + |a_{(2)}||{h_{(1)}}^{(1)}|+ |a_{(2)}||b_{(1)}| \\
  &\quad +| h_{(5)}||g_{(1)}|+ |b_{(1)}||(h_{(2)}| + |b_{(1)}||{h_{(1)}}^{(2)}|+ |{g_{(1)}}^{(1)}| | h_{(4)}|+  |{g_{(1)}}^{(1)}||{h_{(3)}}^{(2)}|
\end{split}\eeq

Therefore,
\begin{multline}\label{alfequalbet}
\alf_0+\alf_1+\alf_2+\alf_3+\alf_4+\alf_5 =  \\
\quad  \bt_0 +\bt_1 +\bt_2+\bt_3-\bt_4-\bt_5-\bt_6+\bt_7-\bt_8+\bt_9-\bt_{10}+\bt_{11}-\bt_{12}+\bt_{13}+\bt_{14}+\bt_{15}.
\end{multline}

By \eqref{yyy}, \eqref{xxx} and \eqref{alfequalbet} we have
\[\Del((a\ot h)\cdot (b\ot g))=\Del(a\ot h)\cdot \Del(b\ot g). \nn \]
The rest of the proof can be done by a similar method.
\epf

\bre
The above lemma is actually true in any symmetric monoidal category. The proof involves braiding diagrams.
\ere

We skip the lengthy, but computational, proof of the following theorem which is the main result of this paper.

\bth
The actions and coactions defined in  Lemmas \eqref{FG2isUg1modalg}  and \eqref{Ug1Fg2comodcoalg} satisfy the conditions \eqref{delhrhda}-\eqref{cpcondbicrsprdct} in Proposition \eqref{supercptbltycondnonsup}.
Therefore we have a bicrossproduct super Hopf algebra $F(G^s_2) \lrbicross U(\mfg^s_1)$ with the following structures:
\[(a\ot h)(b\ot g)= (-1)^{|h_{(2)}||b|} \quad  a(h_{(1)} \rhd b) \ot h_{(2)} g       , \nn \]
\[\Del(a \ot h)=(-1)^{|{h_{(1)}}^{(1)}||{a_{(2)}}| } \quad  a_{(1)} \ot {h_{(1)}}^{(1)} \ot a_{(2)}{h_{(1)}}^{(2)}\ot h_{(2)}       , \nn \]
\[\epn(a \ot h)= \epn(a)\epn(h)     , \nn \]
\[S(a \ot h)= (-1)^{|h^{(1)}||a|} \quad (1\ot S(h^{(1)}))(S(a h^{(2)}) \ot 1)     ,\nn \]
for any $a,b$ in $F(G^s_2)$ and $g,h$ in $U(\mfg^s_1)$.
\ethe

The bicrossproduct super Hopf algebra $\mh^s_1 := F(G^s_2) \lrbicross U(\mfg^s_1)$ is the super version
of the Connes-Moscovici Hopf algebra $\mh_1$.


\bibliographystyle{amsalpha}

\end{document}